\documentclass[11pt]{amsart}
\usepackage{amssymb}
\usepackage{amsfonts}
\usepackage{amsthm}
\numberwithin{equation}{section}
\theoremstyle{plain}
\newtheorem{theorem}{Theorem}[section]

\newtheorem{proposition}[theorem]{Proposition}

\theoremstyle{definition}
\newtheorem{definition}[theorem]{Definition}
\theoremstyle{remark}
\newtheorem*{remark}{Remark}
\newtheorem*{example}{Example}
\newcommand{\refE}[1]{(\ref{E:#1})}
\newcommand{\refS}[1]{Section~\ref{S:#1}}
\newcommand{\refSS}[1]{Section~\ref{SS:#1}}
\newcommand{\refT}[1]{Theorem~\ref{T:#1}}

\newcommand{\refP}[1]{Proposition~\ref{P:#1}}

\newcommand{\R}{\ensuremath{\mathbb{R}}}
\newcommand{\C}{\ensuremath{\mathbb{C}}}
\newcommand{\N}{\ensuremath{\mathbb{N}}}

\renewcommand{\P}{\ensuremath{\mathbb{P}}}
\newcommand{\Z}{\ensuremath{\mathbb{Z}}}

\newcommand{\nl}{\hfill\newline}
\renewcommand{\i}{{\,\mathrm{i}\,}}
\newcommand{\scp}[2]{{\langle #1,#2\rangle}}
\newcommand{\w}{\omega}

\newcommand{\e}{\mathrm{e}}

\newcommand{\ghm}[1][m]{\Gamma_{hol}(M,L^{#1})}
\newcommand{\gh}{\Gamma_{hol}(M,L)}

\newcommand{\gul}{\Gamma_{\infty}(M,L)}
\newcommand{\gulm}[1][m]{\Gamma_{\infty}$(M,L^{#1})}
\renewcommand{\d}{\partial}
\newcommand{\db}{\overline{\partial}}

\newcommand{\zb}{\overline{z}}

\newcommand{\pnc}[1][N]{\ensuremath{\P^{#1}(\C)}}
\newcommand{\volu}{\mathrm{vol}}

\newcommand{\hh}{\hat h}
\def\im{\text{Im\kern1.0pt }}
\def\re{\text{Re\kern1.0pt }}
\newcommand{\Cim}{C^{\infty}(M)}
\newcommand{\cim}{C^{\infty}(M)}
\newcommand{\PGL}{\mathrm{PGL}}
\newcommand{\eam}{e_{\alpha}^{(m)}}
\newcommand{\ebm}{e_{\beta}^{(m)}}
\def\pbar{\overline{\partial}}
\def\la{\lambda}
\def\a{\alpha}

\def\Cim{C^{\infty}(M)}
\def\ghm{\Gamma_{hol}(M,L^{m})}

\def\gh{\Gamma_{hol}(M,L)}

\def\gul{\Gamma_{\infty}(M,L)}
\def\gulm{\Gamma_{\infty}(M,L^{m})}
\def\Lp{{\mathrm{L}}^2(M,L)}
\def\Lpm{{\mathrm{L}}^2(M,L^m)}

\def\Lqv{{\mathrm{L}}^2(Q,\mu)}

\def\Tfm{T_f^{(m)}}
\def\Tgm{T_g^{(m)}}
\def\Tfgm{T_{\{f,g\}}^{(m)}}
\def\Tma#1{T_{#1}^{(m)}}
\def\Hm{\mathcal {H}^{(m)}}
\def\Hc{{\mathcal {H}}}
\def\d{\partial}
\def\db{\overline{\partial}}

\def\Pfz#1{\frac {\partial #1}{\partial z}}
\def\Pfzb#1{\frac {\partial#1}{\partial\overline{z}}}
\def\End{\mathrm{End}}

\def\Tr{\operatorname{Tr}}
\def\volu{\operatorname{vol}}
\def\hh{\hat h}
\renewcommand{\Im}{I^{(m)}}

\newcommand{\eghm}{\mathrm{End}(\Gamma_{hol}(M,L^{(m)}))}
\newcommand{\skp}[2]{{\langle #1,#2\rangle}}

\newcommand{\pimh}{\hat\Pi^{(m)}}
\newcommand{\Bm}{\mathcal{B}_m}
\newcommand{\Pim}{\Pi^{(m)}}
\newcommand{\pim}{\Pi^{(m)}}
\begin{document}
\vspace*{-1cm}
\hspace*{\fill} Mannheimer Manuskripte 257

\hspace*{\fill} math/0009219

\vspace*{2cm}

\title{Berezin-Toeplitz quantization and Berezin transform}
\author{Martin Schlichenmaier}
\address[Martin Schlichenmaier]{Department of Mathematics and 
  Computer Science, University of Mannheim\\ 
       Seminargeb\"aude A5 \\
         D-68131 Mannheim \\
         Germany}
\email{schlichenmaier@math.uni-mannheim.de}
\begin{abstract}
In this lecture results on the Berezin-Toeplitz quantization of 
arbitrary compact quantizable K\"ahler manifolds are presented.
 These results are obtained in
joint work with M.~Bordemann and E.~Meinrenken.
The existence of the Berezin-Toeplitz deformation quantization is
also covered.
Recent results 
obtained in joint work with A. Karabegov on the asymptotic expansion 
of the Berezin transform for arbitrary quantizable compact K\"ahler
manifolds are explained.
As an application the asymptotic expansion of the Fubini-Study
fundamental form under the coherent state embedding is considered.
Some comments on the dynamics of the quantum operators are
given.

(This is an extended write-up of an
invited lecture presented at the workshop
``Asymptotic properties of time evolutions in classical and
quantum systems'', September 13-17, 1999, Bologna, Italy.)
\end{abstract}
\subjclass{}
\keywords{deformation quantization, K\"ahler manifolds, Berezin transform,
semi-classical approximation, Berezin-Toeplitz quantization, geometric
quantization, Toeplitz operators, Bergman kernel}
\date{September 23, 2000}
\maketitle
\section{Introduction}\label{S:intro}
In this lecture I would like to present results on the quantization of compact
K\"ahler mani\-folds.
More precisely, I consider the Berezin-Toeplitz quantization and its 
related Berezin-Toeplitz deformation quantization.
These quantization schemes are very much adapted to the situation 
where the phase-space manifold is a K\"ahler manifold 
$(M,\w)$.
Whereas the basic objects can be defined for arbitrary K\"ahler manifolds
the methods employed in the proofs of the theorems rely on
the compactness of $M$.
This is the reason that I will restrict the presentation 
to compact K\"ahler manifolds from the beginning.

The basic set-up is the set-up of geometric quantization.
Firstly, this means that the K\"ahler form $\w$ endows the algebra
of differentiable functions on $M$ with a Poisson structure.
Secondly, we assume that the K\"ahler manifold is a quantizable
K\"ahler manifold, i.e. it admits a holomorphic 
quantum line bundle with curvature essentially equal to the
K\"ahler form.
To every differentiable function on $M$ there exists a family
of Toeplitz operators $\Tfm$ $(m\in\N_0)$ operating on the space of
global holomorphic sections of the $m$-th tensor power
of the quantum line bundle. 
In this way one obtains the Berezin-Toeplitz quantization.
The set-up is explained in detail in 
\refS{set}.
It was shown by Bordemann, Meinrenken and myself \cite{BMS}
that the Berezin-Toeplitz quantization  has the correct semi-classical
behaviour for $m\to\infty$.
The results are presented in \refS{results}.

The proofs essentially use the theory of generalized Toeplitz 
operators developed by Boutet de Monvel and Guillemin 
\cite{BGTo}.
A sketch of the techniques involved is given in \refS{toeplitz}.
With the same techniques it is possible to construct a
formal deformation quantization (also called a star product)
\cite{Schlhab},\cite{SchlBia95},\cite{Schldef}.
The results are also presented in \refS{results}.

As it is shown in joint work with A. Karabegov \cite{KaraSchl}
the obtained star product has some nice
properties 
(i.e. ``separation of variables'', see below).
On the way to prove these properties it is shown that the
Berezin transform has an asymptotic expansion of certain
type. This results is of independent interest.
The Berezin transform associates to the function
$f$ via its Toeplitz operator $\Tfm$ and the 
covariant Berezin symbol $\sigma(\Tfm)$ of this operator another
differentiable function $I^{(m)}(f)$.
Its asymptotics for $m\to\infty$ carries important information.
Some of the results obtained in \cite{KaraSchl} 
are explained in \refS{bertrans} and in the appendix.

Applications of the
asymptotic expansions are given in 
 \refS{asymp}.
Beside another proof of one part of the correct
semi-classical behaviour of the Berezin-Toeplitz
quantization scheme, the asymptotic expansion of the
pull-back of
the Fubini-Study form of the projective space in which the
manifold is embedded using the global holomorphic sections
of the $m$-th tensor power of the quantum line bundle for
$m\to\infty$ is identified.
The existence of such an asymptotic expansion was 
proven by Zelditch \cite{Zel}.
Here, we show that it can be identified 
(after assigning to the asymptotic expansion 
a formal form) with the 
classifying Karabegov form for the ``opposite dual''
star product to the Berezin-Toeplitz star product.

The results presented  are mainly dealing with kinematics.
It is possible to study dynamics by 
considering operators on the collection of quantum Hilbert
spaces given by the space of global holomorphic sections of
the quantum line bundle and its tensor powers.
Some comments can be found in \refS{dynamics}.
Further work is in progress. There are some results in this
field due to Zelditch \cite{Zelid}, and others.

Finally, I like to thank the organizers of the workshop 
S. Graffi and A. Martinez for the kind invitation 
to present these results, and them and other
participants for stimulating discussions.
I also like to thank M. Engli\v s for bringing the work
of Zelditch to my attention.
\section{The set-up}\label{S:set}
In the following we consider only phase-space manifolds which carry the
structure of a compact
K\"ahler manifold $\ (M,\w)\ $.
In particular, we take as symplectic form the K\"ahler form $\w$ 
which is
a non-degenerate closed positive $(1,1)$-form.
If the complex dimension of $M$ is $n$ then the K\"ahler form $\w$ can be 
written with respect to local holomorphic coordinates
${\{z_i\}}_{i=1,\ldots,n}$ as 
\begin{equation}
\w=\i\sum_{i,j=1}^{n} g_{ij}(z)dz_i\wedge d\bar z_j,
\end{equation}
with local functions $g_{ij}(z)$ such that the 
matrix $(g_{ij}(z))_{i,j=1,\ldots,n}$ is hermitian and
positive definite.

Denote by $\Cim$ the algebra of 
complex-valued (arbitrary often) differentiable functions
with associative product the point-wise multiplication.
Using the K\"ahler form $\w$ 
one  assigns to every  $f\in\Cim$ its {\em Hamiltonian vector field} $X_f$ and 
to every pair of functions $f$ and $g$ the 
{\em Poisson bracket} $\{.,.\}$ via
\begin{equation}
\label{E:Poi}
\w(X_f,\cdot)=df(\cdot),\qquad  
\{\,f,g\,\}:=\w(X_f,X_g)\ .
\end{equation}
With the Poisson bracket 
$\Cim$ 
becomes a {\em Poisson algebra}.
This means  that $\{.,.\}$ defines an additional  Lie algebra structure on
$\Cim$ which fulfills 
as compatibility condition with the associative structure 
the Leibniz rule
$\{fg,h\}=f\{g,h\}+\{f,h\}g$.
 
The K\"ahler manifold $(M,\w)$ is called {\em quantizable}
if there exists an
associated quantum line bundle $(L,h,\nabla)$,
where $L$ is a holomorphic line bundle $L$ over $M$, 
$h$ a Hermitian metric   on $L$, and 
$\nabla$ a connection  compatible with the metric $h$ and the 
complex structure (and hence uniquely fixed),
such that the 
curvature  form of the connection and the K\"ahler form $\w$ of the manifold
are related  as 
\begin{equation}\label{E:quant}
curv_{L,\nabla}(X,Y):=\nabla_X\nabla_Y-\nabla_Y\nabla_X-\nabla_{[X,Y]}
=
-\i\w(X,Y)\ .
\end{equation}
Equation \refE{quant} is called the {\it quantization 
condition}.
If the metric is represented as 
a function $\hh$ 
with respect to 
local complex coordinates  and a 
local holomorphic frame of the bundle
the quantization condition reads as 
$\ \i\db\d\log\hh=\w\ $.

\begin{example}
The Riemann sphere, (i.e. the complex projective line)
 $\P(\C)=\C\cup \{\infty\}\cong S^2$. With respect to the quasi-global
coordinate $z$  the form can be given as
\begin{equation}
\w=\frac {\i}{(1+z\zb)^2}dz\wedge d\zb\ .
\end{equation}
The quantum line bundle  $L$ is the hyperplane  bundle.
For the Poisson bracket one obtains
\begin{equation}
\{f,g\}=\i(1+z\zb)^2\left(\Pfzb f\cdot\Pfz g-\Pfz f\Pfzb g\right)\ .
\end{equation}
\qed
\end{example}

\begin{example}
The (complex-) onedimensional torus $M$.
Up to isomorphy it can be given as
$M = \C/\Gamma_\tau$ where $\ \Gamma_\tau:=\{n+m\tau\mid n,m\in\Z\}$
is a lattice with $\ \im \tau>0$.
As K\"ahler form we take
\begin{equation}
\w=\frac {\i\pi}{\im \tau}dz\wedge d\zb\ ,
\end{equation}
with respect to the coordinate $z$ on the
covering space $\C$.
Clearly this form is invariant under $\Gamma_\tau$ and hence
well-defined on $M$.
The corresponding quantum line bundle is the theta line bundle
of degree 1, i.e. the bundle whose global sections are
scalar multiples of the Riemann theta function.
For the Poisson bracket one obtains
\begin{equation}
\{f,g\}=\i\frac{\im\tau}{\pi}\left(\Pfzb f\cdot\Pfz g-\Pfz f\Pfzb g\right)\ .
\end{equation}
\qed
\end{example}

The quantization condition  \refE{quant} 
implies that $L$ is a positive 
line bundle. By the Kodaira embedding theorem   $L$ 
is  ample, which says that  a certain
tensor power $L^{m_0}$ of $L$ is very ample, i.e. 
the global holomorphic sections of  $L^{m_0}$ 
can be used to embed the phase space manifold $M$ 
into projective
space (the embedding will be explained in a second).
In the following we will assume $L$ already to be  very ample.
If $L$ is not very ample we choose $m_0\in\N$  such that 
the bundle $L^{m_0}$ is very ample and take  
this bundle as quantum line bundle with respect to 
the K\"ahler form $m_0\w$  
on $M$. The   underlying complex manifold structure
is not changed.

We choose  local holomorphic coordinates $z$ for $M$ and
 a local holomorphic frame $e(z)$
for the bundle $L$.
After these choices the basis 
$s_0,s_1,\ldots,s_N$ for
$\gh$, the space of global holomorphic sections of the bundle $L$,
can be uniquely described by
local holomorphic functions
$\hat s_0,\hat s_1,\ldots,\hat s_N$ defined via $s_j(z)=\hat s_j(z)e(z)$.
The embedding is given by the map
\begin{equation}\label{E:embed}
M\hookrightarrow \pnc,\quad z\mapsto \phi(z)=
(\hat s_0(z):\hat s_1(z):\cdots:\hat s_N(z))\ .
\end{equation}
Note that the point $\phi(z)$ 
in projective space neither depends on the choice of
local coordinates nor on the choice of the local frame for the bundle $L$.
A different choice of basis is equivalent unter a $\PGL(N,\C)$ action
on the embedding space.
By this embedding we see that  quantizable compact K\"ahler manifolds are 
as complex manifolds projective
algebraic manifolds.
The converse is also true, see \cite{SchlBia95},\cite{BerSchlcse}.
Note that the embedding is an embedding as complex manifolds
not an isometric embedding as K\"ahler manifolds.

To introduce a scalar product on the space of sections 
we take the Liouville form
$\ \Omega=\frac 1{n!}\w^{\wedge n}\ $ as volume form  on $M$
and set 
for  the scalar product
and the norm
\begin{equation}
\label{E:skp}
\langle\varphi,\psi\rangle:=\int_M h (\varphi,\psi)\;\Omega\  ,
\qquad
||\varphi||:=\sqrt{\langle \varphi,\varphi\rangle}\ ,
\end{equation}
on the space $\gul$ of global $C^\infty$-sections.
Let $\Lp$ be the  L${}^2$-completion of $\gul$,  and
$\gh$ be its (due to the compactness of $M$) 
finite-dimensional closed subspace of global holomorphic
sections.
Let $\ \Pi:\Lp\to\gh\ $ be the projection.
\begin{definition}
For $f\in\Cim$ the {\it Toeplitz operator}  $T_f$
is defined to be
\begin{equation}
T_f:=\Pi\, (f\cdot):\quad\gh\to\gh\ .
\end{equation}
\end{definition}
In words: One takes a holomorphic section $s$ 
and multiplies it with the
differentiable function $f$. 
The resulting section $f\cdot s$ will only be  differentiable.
To obtain a holomorphic section one has to project  it back
on the subspace of holomorphic sections.

The linear map
\begin{equation}
T:\Cim\to \End\big(\gh\big),\qquad  f\to T_f\ ,
\end{equation}
is the  {\it Berezin-Toeplitz quantization map}. 
Because in general
$
T_f\, T_g=\Pi\,(f\cdot)\,\Pi\,(g\cdot)\,\Pi\ne
\Pi\,(fg\cdot)\,\Pi =T_{fg}
$,
it is
neither a Lie algebra homomorphism nor
an associative algebra homomorphism.
The Berezin-Toeplitz quantization
is  a map
from the commutative algebra of functions to a noncommutative
finite-dimensional (matrix) algebra.
The finite-dimensionality is 
due to compactness of $M$.
A lot of classical information will get lost. To recover this
information one should consider not just the bundle $(L,\nabla,h)$ alone but
all its tensor powers $(L^m,\nabla^{(m)},h^{(m)})$
(we use the  notation $L^m:=L^{\otimes m}$)
and apply the above constructions for every $m\in\N_0$.
Note that if
 $\hat h$ corresponds to the metric $h$ with respect to 
 a local holomorphic frame $e$
of the bundle $L$ then $\hat h^m$ corresponds to the metric $h^{(m)}$ 
with respect to the frame $e^{\otimes m}$ for the bundle $L^m$.
In this way one obtains a family of
finite-dimensional(matrix) algebras and a family of maps
\begin{equation}
\Tma {}:\Cim\to \End\big(\ghm\big),\qquad  f\to \Tma f=\Pi^{(m)}(f\cdot)
\ , m\in\N_0\ .
\end{equation}
This infinite family should in some sense ``approximate'' the
algebra $\Cim$.
\section{The semi-classical behaviour and the Berezin-Toeplitz deformation
quantization}
\label{S:results}
For the rest of the lecture let $(M,\omega)$ be a fixed quantizable
compact K\"ahler manifold and $(L,h)$ a fixed very ample quantum line
for it.
Because the connection is uniquely fixed I will drop it in the notation.
Denote for $f\in\Cim$ by $||f||_\infty$ the sup-norm of $\ f\ $ on $M$ and by 
\begin{equation}
||\Tfm||:=\sup_{\substack {s\in\ghm\\ s\ne 0}}\frac {||\Tfm s||}{||s||}
\end{equation}
the operator norm with respect to the analogue of \refE{skp}
on $\ghm$.
The following theorem was shown in 1994.
\begin{theorem}
\label{T:approx}
[Bordemann, Meinrenken, Schlichenmaier]

(a) For every  $\ f\in \Cim\ $ there exists $C>0$ such that   
\begin{equation}\label{E:norma}
||f||_\infty -\frac Cm\quad
\le\quad||\Tfm||\quad\le\quad ||f||_\infty\ .
\end{equation}
In particular, $\lim_{m\to\infty}||\Tfm||= ||f||_\infty$.

(b) For every  $f,g\in \Cim\ $ 
\begin{equation}
\label{E:dirac}
||m\i[\Tfm,\Tgm]-\Tfgm||\quad=\quad O(\frac 1m)\ .
\end{equation}

(c) For every  $f,g\in \Cim\ $ 
\begin{equation}\label{E:prod}
||\Tfm\Tgm-T^{(m)}_{f\cdot g}||\quad=\quad O(\frac 1m)
\ .
\end{equation}
\end{theorem}
These results are contained in  Theorem  4.1, 4.2,
resp. in Section 5 in \cite{BMS}.
Note that (c) 
generalizes trivially to finitely many functions.

In \cite{BHSS}  the notion of $L_{\alpha}$, resp.
$gl(N)$, resp. $su(N)$ quasi-limit was used introduced.
It is related to the concept of continuous fields of $C^*$-algebras
and the notion of strict quantization (see \cite{LandTop},\cite{Riefque})
for the definitions).
It was conjectured in \cite{BHSS} that for every 
compact K\"ahler manifold the Poisson algebra of functions 
is a $gl(N)$ quasi-limit. This was proved in \cite{BMS}.
The result is of special interest in the theory of membranes.
It follows also that the Berezin-Toeplitz quantization is 
a strict quantization.

\begin{remark}
There is another geometric concept of quantization, the {\it geometric
quantization} introduced by Kostant and Souriau.
Due to a result of Tuynman \cite{TuyQ} (see also \cite{BHSS} for a 
coordinate independent proof) 
for compact K\"ahler manifolds 
both quantization schemes  have the same
semi-classical behaviour.
More precisely, if $Q_f^{(m)}$ denotes the well-known  operator
of geometric quantization 
with respect  to K\"ahler polarization, and
$\Delta$ is the Laplacian with respect to the K\"ahler metric given by
$\omega$ then 
$Q_f^{(m)}=\i\cdot T_{f-\frac 1{2m}\Delta f}^{(m)}$.
\end{remark}

\medskip

With the  help of Toeplitz operators it is 
possible to  construct 
a deformation quantization.
A {\em deformation quantization} is given by a {\it star product}.
I will use both terms interchangeable.
To fix the notation and the factors of $\i$ 
let me recall the definition of a star product.
Let $\mathcal{A}=\Cim[[\nu]]$ be the algebra of formal power
 series in the
variable $\nu$ over the algebra $\Cim$. A product $\star$
 on $\mathcal {A}$ is 
called a (formal) star product for $M$ (or for $\Cim$) if it is an
associative $\C[[\nu]]$-linear product which is $\nu$-adically continuous 
such that
\begin{enumerate}
\item
\qquad $\mathcal{ A}/\nu\mathcal {A}\cong\Cim$, i.e.\quad $f\star g
 \bmod \nu=f\cdot g$,
\item
\qquad $\dfrac 1\nu(f\star g-g\star f)\bmod \nu = -\i \{f,g\}$,
\end{enumerate}
where $f,g\in\Cim$.

We can  write 
\begin{equation}
\label{E:cif}
 f\star g=\sum\limits_{j=0}^\infty C_j(f,g)\nu^j\ ,
\end{equation}
with
$ C_j(f,g)\in\Cim$. The $C_j$ should be  $\C$-bilinear in $f$ and $g$.
The conditions 1. and 2.  can 
be reformulated as 
\begin{equation}
\label{E:cifa}
C_0(f,g)=f\cdot g,\qquad\text{and}\qquad
C_1(f,g)-C_1(g,f)=-\i \{f,g\}\ .
\end{equation}
By the $\nu$-adic continuity \refE{cif} fixes $\star$ on  $\mathcal{A}$.
\begin{theorem}
\label{T:star}
There exists a unique (formal) star product $\star_{BT}$ for $M$
\begin{equation}
f \star_{BT} g:=\sum_{j=0}^\infty \nu^j C_j(f,g),\quad C_j(f,g)\in
C^\infty(M),
\end{equation}
in such a way that for  $f,g\in\Cim$ and for every $N\in\N$  we have
with suitable constants $K_N(f,g)$ for all $m$
\begin{equation}
\label{E:sass}
||T_{f}^{(m)}T_{g}^{(m)}-\sum_{0\le j<N}\left(\frac 1m\right)^j
T_{C_j(f,g)}^{(m)}||\le K_N(f,g) \left(\frac 1m\right)^N\ .
\end{equation}
\end{theorem}
This theorem has been proven immediately after 
\cite{BMS} was finished. It has been announced in \cite{SchlBia95},%
\cite{SchlGos}
and the proof was written up in German in  \cite{Schlhab}.
A complete proof published in English 
can be found in \cite{Schldef}.

If we write 
\begin{equation}
\label{E:expans}
T_{f}^{(m)}\cdot T_{g}^{(m)}
\quad \sim\quad \sum_{j=0}^\infty\left(\frac 1m\right)^j
T_{C_j(f,g)}^{(m)}
\qquad (m\to\infty)
\end{equation}
in the following 
we will always assume the strong and precise statement of \refE{sass}.
The same is assumed for other asymptotic formulas appearing 
further down in the article.

In the proof of the above theorems the concept of Toeplitz structure due
to Boutet de Monvel and Guillemin \cite{BGTo} is employed.
In \refS{toeplitz} I will explain this structure and
give a sketch of the proof of Relation \refE{dirac}.
In \refSS{normp} I will give a proof of Relation \refE{norma} 
using results on the Berezin transform. In the original
article \cite{BMS}  a different proof was presented.

A further investigation shows that 
the Berezin-Toeplitz deformation quantization has some important 
properties (see \cite{Schldef} for definitions and proofs).
It is ``null on constants'', i.e. $1\star f=f\star 1=f$.
It is selfadjoint, i.e. $\overline{f\star f}=\overline{g}\star\overline{f}$.
It admits a trace.
As it is shown in \cite{KaraSchl}  it is local and fulfills the 
``separation of variables'' property. Locality means that 
$\mathrm{supp}\, C_j(f,g)\subseteq \mathrm{supp}\, f\cap \mathrm{supp}\, g$
for all $f,g\in\Cim$.
{}From the locality property it follows that 
the $C_j$ are bidifferential operators and that 
the global star product defines
for every open subset $U$ of $M$ a star product for 
the Poisson algebra $C^\infty(U)$.
``Separation of variables'' \cite{Karasep} means that 
 $f\star k=f\cdot k$ and $k\star g=k\cdot g$ for (locally defined)
holomorphic functions $g$, antiholomorphic functions $f$,     
and arbitrary functions $k$.
Note that in Karabegov's notation the r\^oles of the holomorphic
and antiholomorphic functions is switched.
\section{The Toeplitz structure}\label{S:toeplitz}
In \cite{BMS} the set-up for the proof of the 
approximation results was given.
Let me recall for further reference the main definitions.
A  more detailed exposition can be found in 
\cite{Schlhab}.
Take $\ (U,k):=(L^*,h^{-1})\ $ the dual of the very ample quantum line bundle
with dual metric $k$,
$Q:=\{\lambda\in U\mid k(\lambda,\lambda)=1\}$ 
the unit circle bundle inside $U$,   and
$\tau: Q\to M$ the projection.
Note that for the projective space 
with quantum line bundle the hyperplane section bundle  $H$,
the bundle $U$ is just the tautological
bundle. Its  fibre over the point $z\in\P^N(\C)$ consists of
the line in $\C^{N+1}$ which is represented by $z$. In particular,
for the projective space 
the total space of $U$ with  the zero section removed can be identified
with $\C^{N+1}\setminus\{0\}$.
The same picture remains true for the via the very ample quantum
line bundle in projective space embedded manifold $M$.
The quantum line bundle will be the pull-back of $H$ 
(i.e. its restriction to the embedded manifold) and its
dual is the pull-back of the tautological bundle.

In the following we use $E\setminus 0$ to denote the total space of
a vector bundle $E$ with the image of the zero section removed. 
Starting from the function
$\hat k(\la):=k(\la,\la)$ on $U$ we define
 $\tilde a:=\frac {1}{2\i}(\d-\pbar)\log \hat k$ on
$U\setminus 0$ (with respect to the complex
structure on $U$) and denote by $\alpha$
its restriction   to $Q$.
Now $d\a=\tau^*\w$ (with $d=d_Q$) and $\mu=\frac 1{2\pi}\tau^*\Omega
\wedge \a$ is a volume form
on $Q$. 
We have for $f\in\Cim$ the relation $\int_Q(\tau^*f)\mu=\int_Mf\Omega$.
Recall that $\Omega$ is the Liouville volume form on $M$.

With respect to $\mu$ we take the L${}^2$-completion $\Lqv$
of the space of functions on $Q$.
The generalized {\em Hardy space} $\Hc$ is the closure of the
functions in  $\Lqv$ which can be extended to
holomorphic functions on the whole
disc bundle
$D:=\{\lambda\in U\mid k(\lambda,\lambda)\le 1\}$.
The generalized {\em Szeg\"o projector} is the projection
\begin{equation}
\label{E:szproj}
\Pi:\Lqv\to \Hc\ .
\end{equation}
 By the natural circle action the bundle 
$Q$ is a $S^1$-bundle and the tensor powers of $U$ can be
viewed as associated line bundles. The space $\Hc$ is preserved
by the $S^1$-action.
It can be decomposed into eigenspaces 
$\Hc=\prod_{m=0}^\infty \Hm$ where
 $c\in S^1$ acts on $\Hm$ as multiplication
by $c^m$.
The Szeg\"o projector is $S^1$ invariant and 
can be decomposed into its components, the Bergman projectors  
$\pimh:\Lqv\to\Hm$.

Sections of $L^m=U^{-m}$ can be identified with functions $\psi$ on $Q$ which
satisfy the equivariance condition
$\psi(c\la)=c^m\psi(\la)$, i.e. which are homogeneous of degree $m$.
This identification is given via the map
\begin{equation}\label{E:secident}
\gamma_m:\Lpm \to \Lqv,\quad s\mapsto \psi_s\quad\text{where}\quad
\psi_s(\alpha)=\alpha^{\otimes m}(s(\tau(\alpha))),
\end{equation}
which turns out to be an isometry.
Recall that $\Lpm$ has a scalar product given in an corresponding
way to \refE{skp}.
Restricted to the holomorphic sections we obtain the 
isometry 
\begin{equation}\label{E:sechident}
\gamma_m:\ghm \cong \Hm.
\end{equation}. 

There is the notion of Toeplitz structure
$(\Pi,\Sigma)$ as developed by Boutet de Monvel  and
Guillemin
in \cite{BGTo},\cite{GuCT}.
Here 
$\Pi $ is the  Szeg\"o projector  \refE{szproj}
and $\Sigma$  is  the submanifold 
\begin{equation}
\Sigma=\{\;t\alpha(\lambda)\;|\;\lambda\in Q,\,t>0\ \}\ \subset\  T^*Q
\setminus 0
\end{equation}
of the tangent bundle of $Q$ 
defined with the help of the 1-form $\alpha$.
They showed that $\Sigma$ is a 
symplectic submanifold.
A (generalized) {\em Toeplitz operator} of order $k$ is  an operator
$A:\Hc\to\Hc$ of the form
$\ A=\Pi\cdot R\cdot \Pi\ $ where $R$ is a
pseudodifferential operator ($\Psi$DO)
of order $k$ on
$Q$.
The Toeplitz operators build a ring.
The  symbol of $A$ is the restriction of the
principal symbol of $R$ (which lives on $T^*Q$) to $\Sigma$.
Note that $R$ is not fixed by $A$, but
Guillemin and Boutet de Monvel showed that the  symbols
are well-defined and that they obey the same rules as the
symbols of   $\Psi$DOs.
In particular we have the following relations
\begin{equation}
\label{E:symbol}
\sigma(A_1A_2)=\sigma(A_1)\sigma(A_2),\qquad
\sigma([A_1,A_2])=\i\{\sigma(A_1),\sigma(A_2)\}_\Sigma.
\end{equation}
Here $\{.,.\}_{\Sigma}$ is the restriction of the canonical
Poisson structure of $T^*Q$ to $\Sigma$ coming from the
canonical symplectic form $\w_0$ on $T^*Q$.

We have to deal with 
two Toeplitz operators:
\nl
(1) The generator of the circle action
gives the  operator $D_\varphi=\dfrac 1{\i}\dfrac {\partial}
{\partial\varphi}$. It is an operator of order 1 with symbol $t$.
It operates on $\Hm$ as multiplication by $m$.
\nl
(2) For $f\in\Cim$ let $M_f$ be the  operator on
$\Lqv$ 
corresponding to multiplication with $\tau^*f$.
We set%
\footnote{
There should be no confusion with the operator
$T_f=T_f^{(1)}$ introduced above.}
$\ T_f=\Pi\cdot M_f\cdot\Pi:\Hc\to\Hc\ $.
Because $M_f$ is constant along the fibres of $\tau$, $T_f$ 
commutes with the circle action.
Hence
$\ T_f=\prod\limits_{m=0}^\infty\Tfm\ $,
where $\Tfm$ denotes the restriction of $T_f$ to $\Hm$.
After the identification of $\Hm$ with $\ghm$ we see that these $\Tfm$
are exactly the Toeplitz operators  $\Tfm$ introduced in \refS{set}.
In this sense  $T_f$ is called  the global Toeplitz operator and
the $\Tfm$ the local Toeplitz operators.
$T_f$ is an operator of order $0$.
Let us denote by
$\ \tau_\Sigma:\Sigma\subseteq T^*Q\to Q\to M$ the composition
then we obtain for the symbol  $\sigma(T_f)=\tau^*_\Sigma(f)$.

Now we are able to proof \refE{dirac}.
The commutator
$[T_f,T_g]$ is a  Toeplitz operator of order $-1$.
Using $\ {\omega_0}_{|t\alpha(\lambda)}=-t\tau_\Sigma^*\omega\ $
if $t$ is a fixed positive number, we obtain
with \refE{symbol} that the 
principal symbol of the commutator equals
\begin{equation}
\sigma([T_f,T_g])(t\alpha(\lambda))=\i\{\tau_\Sigma^* f,\tau_\Sigma^*g
\}_\Sigma(t\alpha(\lambda))=
-\i t^{-1}\{f,g\}_M(\tau(\lambda))\ .
\end{equation}
Now consider the Toeplitz operator
\begin{equation}
A:=D_\varphi^2\,[T_f,T_g]+\i D_\varphi\, T_{\{f,g\}}\ .
\end{equation}
Formally this is an operator of order 1.
Using $\ \sigma(T_{\{f,g\}})=\tau^*_\Sigma \{f,g\}$ 
and $\sigma(D_\varphi)=t$ we see that its principal
symbol vanishes. Hence
it is an operator of order 0.
Now $M$ and hence $Q$ are compact manifolds.
 This implies that $A$ is a bounded
operator ($\Psi$DOs of order 0 on compact manifolds are bounded).
It is obviously $S^1$-invariant and we can write
$A=\prod_{m=0}^\infty A^{(m)}$
where $A^{(m)}$ is the restriction of $A$ on the space $\Hm$.
For the norms we get $\ ||A^{(m)}||\le ||A||$.
But
\begin{equation}
A^{(m)}=A_{|\Hm}=m^2[\Tfm,\Tgm]+\i m\Tfgm.
\end{equation}
Taking the norm bound and dividing it by $m$ we get the 
the statement \refE{dirac}.

Quite similar \refE{prod} and \refT{star} can be proved.
The original proof of \refE{norma} presented in \cite{BMS} uses
different techniques. Below I will introduce the Berezin transform.
With the help of its asymptotic expansion shown in 
\cite{KaraSchl} I will give further down a different proof of 
\refE{norma}.
\section{Berezin symbols and the Berezin transform}\label{S:bertrans}
\subsection{Coherent States}\label{SS:coherent}
Let the situation as in the previous section. In particular $L$ is assumed
 to be
already very ample, $U=L^*$ is the dual of the quantum line bundle, 
$Q\subset U$ the unit circle bundle, and
$\tau:Q\to M$ the projection.
Recall that the sections of 
$L^m$ can be identified with the homogeneous functions on 
$Q$  of degree $m$, see \refE{secident}.
For every element $\alpha\in U\setminus 0$ there exists a unique section
$\eam\in\ghm$ such that
\begin{equation}
\skp{s}{\eam}=\psi_s(\alpha)=\alpha^{\otimes m}(s(\tau(\alpha)))
\end{equation}
for all $s\in\ghm$.
This section is called the coherent vector associated to the point $\alpha$.
Recall that $\skp{.}{.}$ denotes the scalar product on the space of
global sections $\gulm$.
We assume it to be linear in the first argument and anti-linear in
the second argument.
The definition is dual to the definition of 
coherent vectors of 
Berezin in its
coordinate independent version and extension 
due to Rawnsley \cite{Raw},\cite{CGR1},
see also \cite{BHSS},\cite{Schlhab}.
There the coherent vectors are parameterized by the elements of
$L\setminus 0$.
Our definition has the advantage that we can consider
all tensor powers of $L$ together.

The coherent vectors are antiholomorphic in $\alpha$ and fulfill
\begin{equation}\label{E:cohtrans}
e_{c\alpha}^{(m)}=\bar c\cdot \eam,\qquad c\in\C^*:=\C\setminus\{0\}\ .
\end{equation}
Note that $\eam\equiv 0$ would imply 
that all sections will vanish at the point $x=\tau(\alpha)$.
But this is a contradiction to the very-ampleness of $L$.
Hence, $\eam\not\equiv 0$ and due to \refE{cohtrans} the element 
$[\eam]:=\{s\in\ghm\mid \exists c\in\C^*:s=c\cdot \eam\}$
is a well-defined
element of the projective space $\P(\ghm)$ only depending on 
$x=\tau(\alpha)\in M$.
It is called the {\em coherent state associated to $x\in M$}.

The {\em coherent state embedding} is the 
antiholomorphic embedding
\begin{equation}\label{E:cohemb}
M\quad \to\quad \P(\ghm)\ \cong\ \pnc[N],
\qquad 
x\mapsto [e_{\tau^{-1}(x)}].
\end{equation}
In abuse of notation in this context we will understand
under $\tau^{-1}(x)$ always a non-zero element of the 
fiber over $x$.
The coherent state embedding is up to conjugation the  embedding
with respect to  an orthonormal basis of the sections.
See \cite{BerSchlcse}  for further considerations of the geometry involved.

\subsection{Berezin symbols}\label{SS:symbols}
The {\em covariant Berezin symbol} $\sigma(A)$
of an operator $A\in\eghm$ is defined as
\begin{equation}\label{E:covB}
\sigma(A):M\to\C,\qquad
x\mapsto \sigma(A)(x):=
\frac {\skp {A\eam}{\eam}}{\skp {\eam}{\eam}},\quad \alpha\in\tau^{-1}(x),
\alpha\ne 0\ .
\end{equation}
It is a well-defined function on $M$.

It is shown in \cite{SchlBia98} that the Toeplitz map
$f\mapsto \Tfm$ and the symbol map $A\mapsto \sigma(A)$ are adjoint if one
takes for the operators the Hilbert-Schmidt norm and for the functions 
the  Liouville
measure modified by Rawnsley's epsilon function.
 Here I will not go into the details of the correspondence.
Let me only point out that every operator of $\ghm$ can be represented
as $\Tfm$ with a suitable function $f$ which is in general not
unique.
This $f$ is also called a contravariant symbol $\check \sigma (A)$
of the operator $A$.

Starting from $f\in\cim$ we can assign to it its Toeplitz operator
$\Tfm\in\eghm$ and then assign to $\Tfm$ the covariant symbol 
$\sigma(\Tfm)$. It 
is again an element of $\cim$.  
Altogether we obtain a map $f\mapsto \Im(f):=\sigma(\Tfm)$.
\begin{definition}
The map 
\begin{equation}
\Cim\to\Cim,\qquad
f\mapsto
\Im(f):=\sigma(\Tfm)
\end{equation}
is called Berezin transform.
\end{definition}

{}From the point of view of Berezin's approach \cite{Berebd}
the operator $T_f^{(m)}$ has as a  contravariant symbol $f$.
Hence $\Im$ gives a correspondence between contravariant symbols
and covariant symbols of operators.
The Berezin transform was introduced and studied by 
Berezin \cite{Berebd} for certain classical symmetric 
domains in $\C^n$. These results where 
extended by Unterberger and Upmeier \cite{UnUp},
see also Engli\v s \cite{Eng},\cite{Eng2},\cite{Engbk} 
and Engli\v s and Peetre \cite{EnPe}.
Obviously, the Berezin transform makes also sense
in the compact K\"ahler case which we
consider here.

\subsection{Asymptotic expansion of the Berezin transform}\label{SS:asymp}
The results presented in this subsection are
joint work with Alexander Karabegov \cite{KaraSchl}.
Recall from \refS{toeplitz}
the Szeg\"o projectors 
$ \Pi:\Lqv\to \Hc$ and its components 
$\pimh:\Lqv\to\Hm$, the Bergman projectors.
The Bergman projectors have a smooth integral kernels,
the Bergman kernels 
$\Bm(\alpha,\beta)$ defined on $Q\times Q$, i.e.
\begin{equation}
\pimh(\psi)(\alpha)=\int_Q\Bm(\alpha,\beta)\psi(\beta)\mu(\beta).
\end{equation}
The Bergman kernels can be expressed with the help of the
coherent vectors.
\begin{proposition}\label{P:kernel}
\begin{equation}
\Bm(\alpha,\beta)=\psi_{\ebm}(\alpha)=
\overline{\psi_{\eam}(\beta)}=\skp{\ebm}{\eam}.
\end{equation}
\end{proposition}
\begin{proof}
First note that for $s\in\Lpm$ we have 
$\pimh\gamma_m(s)=\gamma_m\pim s$, i.e.
$\pimh\psi_s=\psi_{\pim s}$.
Due to the fact that $\eam$ is a holomorphic section
\begin{equation*} 
\skp{s}{\eam}=\skp{s}{\pim\eam}=\skp{\pim s}{\eam}
=\psi_{\pim s}(\alpha)=
\pimh\psi_s(\alpha)
=\int_Q\Bm(\alpha,\beta)\psi_s(\beta)\mu(\beta)\ .
\end{equation*}
By the isometry
\begin{equation*}
\skp{s}{\eam}=\skp{\psi_s}{\psi_{\eam}}=
\int_Q\psi_s(\beta)\overline{\psi_{\eam}(\beta)}\mu(\beta)\ .
\end{equation*}
If we compare this two expressions and take the
definition of the coherent vectors we obtain
\begin{equation*}
\Bm(\alpha,\beta)=\overline{\psi_{\eam}(\beta)}=
\overline{\skp{\eam}{\ebm}}=\skp{\ebm}{\eam}
=\psi_{\ebm}(\alpha).
\end{equation*}
\end{proof}
Let $x,y\in M$ and choose $\alpha,\beta\in Q$ with
$\tau(\alpha)=x$ and $\tau(\beta)=y$ then the functions
\begin{equation}
u_m(x):=\Bm(\alpha,\alpha)=
\skp{\eam}{\eam},
\end{equation}
\begin{equation}
v_m(x,y):=\Bm(\alpha,\beta)\cdot \Bm(\beta,\alpha)=
\skp{\ebm}{\eam}\cdot \skp{\eam}{\ebm}
\end{equation}
are well-defined on $M$.
The following proposition gives an integral representation of the
Berezin transform.
\begin{proposition}\label{P:kernelint}
\begin{equation}\label{E:kernelint}
\begin{aligned}
\left(\Im(f)\right)(x)&=\frac 1{\Bm(\alpha,\alpha)}
\int_Q \Bm(\alpha,\beta)\Bm(\beta,\alpha)\tau^*f(\beta)\mu(\beta)
\\
&=
\frac 1{u_m(x)}
\int_M v_m(x,y)f(y)\Omega(y)\ .
\end{aligned}
\end{equation}
\end{proposition}
\begin{proof}
Take any $\alpha\in\tau^{-1}(x)$ with $\alpha\in Q$.
Denote by $M_f$ the operator of pointwise multiplication
of the sections with the function $f$.
\begin{equation}
\begin{aligned}
\bigl(I^{(m)}
f\bigr)(x)&=\sigma\bigl(\Tfm\bigr)(x)=\frac{\skp{\Tfm\eam}{\eam}}{
\scp{\eam}{\eam}}\\
&=\frac{\skp{\Pim M_f\Pim \eam}{\eam}}{\Bm(\alpha,\alpha)}=
\frac{\skp{M_f\eam}{\eam}}{\Bm(\alpha,\alpha)}
\end{aligned}
\end{equation}
Using the isometry \refE{secident} and \refP{kernel}
 we can rewrite the last expression
and obtain
\begin{equation}
\begin{aligned}
\bigl(I^{(m)}
f\bigr)(x)&=
\frac{\skp{(\tau^*f)\psi_{\eam}}{\psi_{\eam}}}{\Bm(\alpha,\alpha)}= 
\frac{1}{\Bm(\alpha,\alpha)}\int_Q (\tau^*f) 
\psi_{\eam}(\beta)\overline{\psi_{\eam}(\beta)}\mu(\beta)\\ 
&=\frac{1}{\Bm(\alpha,\alpha)}\int_Q\Bm(\alpha,\beta) 
\Bm(\beta,\alpha)(\tau^*f)(\beta)\mu(\beta)
=\frac{1}{u_m(x)}\int_M v_m(x,y) f(y)\Omega(y).
\end{aligned}
\end{equation}
\end{proof}
In \cite{KaraSchl} based on works of Boutet de Monvel and Sj\"ostrand 
\cite{BS} on the Szeg\"o kernel and in generalization of a result of
Zelditch \cite{Zel} on the Bergman kernel 
the integral representation 
is used to prove
the existence and the
form of the 
asymptotic expansion of the Berezin transform.
It can be identified with a
formal Berezin transform introduced by A. Karabegov.
Karabegov developed a theory for such formal 
deformation quantizations of (pseudo-) K\"ahler manifolds
which fulfill the  ``separation of variables'' properties. 
In particular, he assigns to every such deformation quantization a
{\em formal Berezin transform} $I$. It is 
a formal power series in the variable
$\nu$ and can be written as
$I=\sum_{j=0}^\infty I_j\nu^j$ with  operators 
$I_j:\Cim\to\Cim$. He proved that it starts with
$I=id+\nu\Delta+\nu^2\cdots$.
In \cite{KaraSchl} it is shown that if we replace 
$\frac 1m$ by the formal variable $\nu$ in the asymptotic  expansion of
the 
Berezin transform $\Im f(x)$ we obtain
the formal Berezin transform  $I(f)(x)$ with respect to a specified 
star product (Theorem 5.9 in \cite{KaraSchl}).
In particular,  we obtain for the asymptotic expansion 
of the Berezin transform at a fixed point $x\in M$ 
\begin{equation}\label{E:Ima}
\Im f(x)\sim f(x)+\frac 1m\Delta f(x)+\cdots,\quad 
\text{for\ }m\to\infty.
\end{equation}
Some details and intermediate results, e.g. a stationary
phase integral \refE{expansion} which is asymptotically equivalent to
$\bigl(I^{(m)}f\bigr)(x)$, are given in an appendix to this 
write-up.
\section{Applications of the asymptotics of the Berezin transform}
\label{S:asymp}
\subsection{Norm preservation}\label{SS:normp}
Here I like to prove that Relation \refE{norma} in 
\refT{approx} can be easily deduced from \refE{Ima}.
First note that 
\begin{equation}\label{E:symine}
||\Im(f)||_\infty=||\sigma(\Tfm)||_\infty\quad\le \quad||\Tfm||\quad\le\quad   
||f||_\infty\ .
\end{equation}
This has been shown in \cite{SchlBia98}.
For the convenience of the reader I will repeat the proof here.
Using Cauchy-Schwarz inequality we calculate ($x=\tau(\alpha)$)
\begin{equation}
| \sigma(\Tfm)(x)|^2=
\frac {|\skp {\Tfm\eam}{\eam}|^2}{{\skp {\eam}{\eam}}^2}\le
\frac {\skp {\Tfm\eam}{\Tfm\eam}}{\skp {\eam}{\eam}}\le
||\Tfm||^2\ .
\end{equation}
Here the last inequality  follows from the definition of the operator norm.
This shows the first inequality in \refE{symine}.
For the second inequality introduce the multiplication
operator $M_f^{(m)}$ on $\gulm$. Then 
$\ ||\Tfm||=||\Pim\,M_f^{(m)}\,\Pim||\le ||M_f^{(m)}||\ $ and
for  $\varphi\in\gulm$,  $\varphi\ne 0$
\begin{equation}
\frac {{||M_f^{(m)} \varphi||}^2}{||\varphi||^2}=
\frac {\int_M h^{(m)}(f \varphi,f \varphi)\Omega}
 {\int_M h^{(m)}(\varphi,\varphi)\Omega}
=
\frac {\int_M f(z)\overline{f(z)}h^{(m)}(\varphi,\varphi)\Omega}
{\int_M h^{(m)}(\varphi,\varphi)\Omega}
\le
||f||{}_\infty^2\ .  
\end{equation}
Hence,
\begin{equation}
||\Tfm||\le ||M_f^{(m)}||=\sup_{\varphi\ne 0}
\frac {||M_f^{(m)}\varphi||}{||\varphi||}\le ||f||_\infty .
\end{equation}

Now choose as $x_e\in M$ a point with $|f(x_e)|=||f||_\infty$.
From the fact that the formal Berezin
transform has
as leading term the identity   it follows that
$\ |(I^{(m)}f)(x_e)-f(x_e)|\le C/m\ $ with a suitable constant
$C$.
This implies
$\ \left| |f(x_e)|-|(I^{(m)}f)(x_e)| \right| \le C/m\ $
and hence
\begin{equation}\label{E:absch}
||f||_\infty-\frac Cm=|f(x_e)|-\frac Cm\quad\le\quad
|(I^{(m)}f)(x_e)|\quad\le\quad ||I^{(m)}f||_\infty\ .
\end{equation}
Putting \refE{symine} and \refE{absch} together we obtain
\begin{equation}\label{E:thma}
||f||_\infty-\frac Cm\quad\le\quad ||T_f^{(m)}||\quad\le\quad
|f|_\infty\ .
\end{equation}
Note that we obtain in this way 
another proof of \cite{BMS}, Theorem 4.1.
\subsection{Pullback of the Fubini-Study form}\label{SS:pfsm}
Starting from the  K\"ahler manifold $(M,\omega)$ and after choosing an
orthonormal basis of the space $\ghm$ we obtain an embedding 
$\phi^{(m)}:M\to\P^{N(m)}$ 
of $M$ into projective space of dimension $N(m)$.
On  $\P^{N(m)}$ we have the standard K\"ahler form, the Fubini-Study  form
$\omega_{FS}$. 
The pull-back
$(\phi^{(m)})^*\omega_{FS}$ will not depend on the orthogonal basis
chosen for the embedding but
in general it will not coincide with 
a scalar multiple of the K\"ahler form $\omega$ we started with
(see \cite{BerSchlcse} for a thorough discussion of the
situation).
It was shown by Zelditch \cite{Zel} by generalizing a result
of Tian that  
$(\Phi^{(m)})^*\omega_{FS}$ admits a complete asymptotic expansion in 
powers of $\frac 1m$ as $m\to\infty$.

The Karabegov classification of star products 
with the ``separation of variables'' property
for $(M,\omega)$ can
be given by assigning to it a unique
formal form 
$\widehat{\omega}=\frac 1\nu\omega+\omega_0+\nu\omega_1+\cdots$
where the  $\omega_i$ for $i\in\N_0$ are closed $(1,1)$-forms on $M$.
From the proof of the identification of the Berezin transform with the 
formal Berezin transform assigned to a certain star product $\star$, 
it follows that the classifying form of 
of exactly this $\star$ coincides with 
the form obtained via the asymptotic expansion  of 
$(\phi^{(m)})^*\omega_{FS}$ if one replaces $\frac 1m$ by $\nu$
(see also the appendix).
Let me add that in Karabegov's theory the ``opposite of the dual star
product''
to $\star$ is the Berezin-Toeplitz star product $\star_{BT}$ of 
\refT{star}.
\section{Some Comments on the dynamics}\label{S:dynamics}
First note that we have 
the following proposition.
\begin{proposition}
The $\Tfm\in\eghm$ fulfill the relation
\begin{equation}
\label{E:Tconj}
{T_f^{(m)}}^*=T_{\overline{f}}^{(m)}\ .
\end{equation}
\end{proposition}
\begin{proof}
Let  $s,t\in \ghm$.
For the scalar product we calculate
\begin{equation}
\skp {s}{T_f^{(m)}t}=
\skp {s}{\Pi^{(m)}ft}=
\skp {s}{ft}=
\skp {\overline{f}s}{t}=
\skp {T_{\overline{f}}^{(m)}s}{t}\ .
\end{equation}
\end{proof}
Due to the fact that the identification of $\ghm$ with $\Hm$ is an isometry,
the Equation \refE{Tconj} is also true for the components of the global
Toeplitz operator. In particular, we have 
${T_f}^*=T_{\overline{f}}$.
This implies that for real-valued functions $f$ on $M$ (e.g.
for a Hamiltonian) the
quantum  operators $T_f$ and $\Tfm$ are self-adjoint operators.

From the theory of generalized Toeplitz operators 
follows some consequences for the spectral
asymptotics of these operators.
Set $d(m):=\dim\Hm$    and let 
$\lambda_1^{(m)},\lambda_2^{(m)},\ldots,\lambda_{d(m)}^{(m)}$ 
be the eigenvalues of the restriction  of 
$T_f$ on $\Hm$.
In particular, these are also the eigenvalues 
of $\Tfm$ on
$\ghm$.
Following \cite{BGTo}  ($n=\dim_{\C} M$)
let  
\begin{equation}
\rho_m:=\frac 1{m^n}\sum_{i=1}^{d(m)}\delta(\lambda-\lambda_i^{(m)})
\end{equation}
be the discrete spectral measure.
By Theorem~13.13 of \cite{BGTo}
it converges weakly to the limit
measure 
\begin{equation}
\rho(g)=\gamma_M\int_M g(f(z))\,\Omega(z),
\end{equation}
with an universal constant $\gamma_M$ only
depending on the manifold $M$.
For $g\equiv 1$ we obtain
\begin{equation}
\label{E:spectr}
\frac 1{m^n}\sum_{i=1}^{d(m)}\lambda_i^{(m)}
=\frac 1{m^n}\Tr^{(m)} \Tfm=
\gamma_M\int_M f\,\Omega+O(\frac 1m)\ .
\end{equation} 
Here $\Tr^{(m)}$ denotes the 
trace on $\End(\ghm)$.
The constant evaluates to 
$\gamma_M={\volu(\P^{n}(\C))}^{-1}$.
By linearity this  extends to complex-valued functions. Hence,
\begin{proposition} \cite{BMS}
Let $f\in\Cim$ and let $n=\dim_{\C}M$. Then 
\begin{equation}
\Tr^{(m)}\,(T^{(m)}_{f})
        =m^n\left(\frac 1{\volu (\P^{n}(\C))}
\int_M f\, \Omega +O(m^{-1})\right)\ .
\end{equation}
\end{proposition}
Indeed a closer analysis (i.e. the application of Equation 13.13 in 
\cite{BGTo}) shows that $\Tr^{(m)}\,(T^{(m)}_{f})$ admits a
complete asymptotic expansion which allows to
construct a formal trace for the Berezin-Toeplitz star 
product \cite{Schldef}.

Of course, now the question arrises how to quantize 
symplectic maps $\Psi$ of the
phase-space (K\"ahler) manifold $(M,\w)$.
A first condition is that the map lifts to 
a contact transformation $\hat \Psi$ of the $S^1$-bundle
$(Q,\alpha)$. Note that $d\alpha=\tau^*\w$ and $(Q,\alpha)$ 
is a contact structure.
Let $L_{\hat \Psi}$ be the translation operator by $\Psi$ then
the Toeplitz operator $\Pi L_{\hat \Psi} \Pi$ is not necessarily a 
unitary operator. But it was shown by Zelditch \cite{Zelid} that   
there exists always a function $f\in \Cim$, such that 
\begin{equation}
U_{\hat \Psi}=
\Pi M_fL_{\hat \Psi} \Pi
\end{equation}
is unitary and commutes with the $S^1$ action.
In particular, it decomposes again as 
\begin{equation}
U_{\hat \Psi}=\prod_{i=0}^\infty
U_{\hat \Psi}^{(m)},
\end{equation}
where each $U_{\hat \Psi}^{(m)}$ is an unitary operator on 
$\Hm\cong\ghm$.
Such kind of maps where studied by Zelditch \cite{Zelid}.
Further work is in progress.
\appendix
\section{More details on the asymptotic expansion}
In this appendix I will express the integral  formula \refE{kernelint}
of the Berezin transform 
up to asymptotic equivalence as a stationary phase integral
\refE{expansion}.
Clearly, here the steps can only be sketched.
For a complete derivation, see \cite{KaraSchl}.
We start with the formula \refE{kernelint} and fix a point $x\in M$. 
Let $W$ be a small contractible 
neighbourhood of $x$. Split $M$ into the two subsets
$W$ and $M\setminus W$ and correspondingly the 
integral. From results of Boutet de Monvel and Sj\"ostrand \cite{BS}
it follows that for $y\in M\setminus W$ the Bergman kernel 
$\Bm(x,y)$ is rapidly
decreasing for $m\to\infty$. For the Bergman kernel on the 
diagonal 
Zelditch \cite{Zel} proved that $u_m(x)=\Bm(x,x)$ expands in  an 
asymptotic series
\begin{equation}\label{E:um}
u_m(x)\sim m^n\sum_{r\ge 0}\left(\frac 1m\right)^rb_r(x), \quad\text{with\ }
b_0=1,\quad\text{for\ } m\to \infty.
\end{equation} 
Hence, the integral over $M\setminus W$ will be rapidly decreasing
as $m\to\infty$ (here the compactness of $M$ is used).
It follows that up to asymptotic equivalence it is enough to 
consider the integral over $W$.
In the following I will use almost analytic extensions. 
Without too much misconception the 
reader 
might imagine 
almost analytic extensions similar to analytic extensions
of a real-analytic object
$a(x)$ defined on $W$ to an object $\tilde a(x,y)$ on $W\times W$ which
is holomorphic in the first variable, antiholomorphic in
the second, and fulfills $\tilde a(x,x)=a(x)$ on $W$.
Such an almost analytic expansion exists for every $C^\infty$-function.

Let $\Phi$ be the (local) K\"ahler potential of $\w$, i.e.
$\w=-\i\partial\bar\partial \Phi$, over 
$W$ (shrinking $W$ if necessary) and $\tilde\Phi$
an almost analytic extension fulfilling 
 $\tilde\Phi(y,x)=\overline{\tilde\Phi(x,y)}$.
Set
\begin{equation}
\chi(x,y):=\tilde\Phi(x,y)-\frac 12\Phi(x)-
-\frac 12\Phi(y)\quad\text{and}\quad
D(x,y):=\chi(x,y)+\chi(y,x)\ .
\end{equation}
Let $e$ be a local holomorphic frame of $U$ over $W$. For any $x\in W$ take
\begin{equation}
\alpha(x):=\frac {e(x)}{\sqrt{k(e(x),e(x))}}
\end{equation}
as corresponding point in $Q_{|W}$.
Now by Theorem 5.6 of \cite{KaraSchl} we obtain that 
there exists an asymptotic expansion of the 
Bergman kernel over $W\times W$ for   $m\to\infty$ 
\begin{equation}\label{E:bergman}
\Bm(\alpha(x),\alpha(y))
\sim m^n\e^{m\chi(x,y)}\sum_{r\ge 0}\left(\frac 1m\right)^r
\tilde b_r(x,y),
\end{equation}
where the 
$\tilde b_r(x,y)$ are almost analytic extensions of
the $b_r(x)$ which are given by \refE{um}. 
Let $b(x,y,m)\in S^0(W\times W\times \R)$ be a 
symbol such that it has the asymptotic expansion
\begin{equation}
b(x,y,m)\sim\sum_{r\ge 0} \left(\frac 1m\right)^r\tilde b_r(x,y).
\end{equation}
We conclude from \refE{bergman} that
we have the asymptotic expansion 
of the two-point function
\begin{equation}
v_m(x,y)\sim m^{2n} e^{mD(x,y)}b(x,y,m)b(y,x,m)
\end{equation}

Now $\phi^x(y):=D(x,y)$ are phase functions such
that $y=x$ are nondegenerate isolated critical points.
Again by shrinking $W$ if necessary, we can achieve that
$y=x$ are the only critical points and that
$\tilde b_0(x,y)$ does not vanish on the closure of $W\times W$.
We obtain the following asymptotic equivalence
\begin{equation}\label{E:expansion}
(I^{m}f)(x)\sim
m^n\int_We^{m\phi^x(y)}f(y)\frac {b(x,y,m)b(y,x,m)}{b(x,x,m)}
\Omega(y)\ .
\end{equation}
This is a stationary phase integral.
From this follows the asymptotic expansion of the
Berezin transform and 
after some additional work \cite{KaraSchl} its identification with the formal
Berezin transform associated to a certain formal star product $\star$.

Let me indicate some steps necessary to identify 
the star product $\star$. If we replace in $\sum_{r\ge 0}\left(\frac 1m\right)^r b_r$
from \refE{um} 
the expression $\frac 1m$ by the formal variable $\nu$ we can find another
formal function $s$ such that $\sum_{r\ge 0}\nu^rb_r(x)=\e^{s(x)}$.
Note that $b_0(x)\equiv 1$.
Now we set $\widehat{\Phi}:=\frac 1\nu\Phi+s$
(again a formal power series).
It turns out that the Karabegov form  $\tilde \omega$ classifying
the star product $\star$ is given by
$
\widehat{\w}=-\i\partial\bar\partial\widehat{\Phi}
$.
The form $\tilde\w$ is a formal Laurent series in the variable $\nu$
\begin{equation}\label{E:wt}
\widehat{\w}=\frac 1\nu\w+\sum_{j=0}^\infty \nu^j\w_j, 
\end{equation}
where $\w$ is the K\"ahler form we started with and the $\w_i$ are
global closed $(1,1)$-forms (not necessarily nondegenerate).

As is shown in \cite{KaraSchl} the Karabegov 
classifying formal form  $\w'_{BT}$ of the  opposite of the 
Berezin-Toeplitz star product $\star_{BT}$ is
\begin{equation}
\w'_{BT}=-\frac 1\nu\w+\w_{can},
\end{equation}
where $\w_{can}$ is the curvature form 
of the canonical holomorphic line bundle of $M$ with fibre metric
given by the metric coming from the Liouville form.
The opposite star product $\star'_{BT}$ is defined as
$f\star'_{BT} g:=g\star_{BT}f$. This switches the r\^oles
played by the holomorphic and antiholomorphic functions in the 
definition of "separation of variables" in accordance with
its use by Karabegov \cite{Karasep}.

Finally, we conclude 
\begin{equation}
cl(\star_{BT})=\frac {1}{\i}\left(\frac 1\nu[\w]-\frac {\epsilon}{2}\right)
\end{equation} 
for the characteristic class of the star product $\star_{BT}$, where
$\epsilon$ is the canonical class of the manifold $M$.

Let me indicate how the result on the asymptotic of
the pull-back of Fubini-Study form 
presented in \refSS{pfsm} follows from the above presented steps.
The pull-back can be given as
\cite[Prop.9]{Zel}
\begin{equation}
\left(\phi^{(m)}\right)^*\w_{FS}=m\w+\i\partial\bar\partial\log u_m(x)\ .
\end{equation}
By the asymptotic expansion \refE{um} of $u_m(x)$ in terms of the 
$b_r$ due to Zelditch, and 
the definition of the formal K\"ahler potential $\widehat{\Phi}$ 
for the classifying form $\widehat{\w}$ via the formal function
$s$ which is defined
above via the $b_r$,
we  conclude that $\widehat{\w}$ can be obtained from
the asymptotic expansion of  
$\left(\phi^{(m)}\right)^*\w_{FS}$ if one replaces $\frac 1m $ by the
formal variable $\nu$.

\providecommand{\bysame}{\leavevmode\hbox to3em{\hrulefill}\thinspace}

\end{document}